\documentclass[reqno]{amsart}

\usepackage{amsfonts,amssymb,graphicx,float}
\usepackage[arrow,matrix,curve]{xy}

\restylefloat{figure}

\makeatletter                                                  %
\def\th@plain{\slshape}                                        %
\makeatother                                                   %

\newcommand{\oi}{[0,1]}
\newcommand{\ooii}{\{0,1\}}
\newcommand{\Nbb}{\mathbb{N}}
\newcommand{\Zbb}{\mathbb{Z}}

\newcommand{\Rbb}{\mathbb{R}}
\newcommand{\one}{{\rm 1\mskip-4mu l}}

\newcommand{\monk}{M\"onkemeyer}

\newcommand{\Bcal}{\mathcal{B}}
\newcommand{\Fcal}{\mathcal{F}}

\newcommand{\Mcal}{\mathcal{M}}
\newcommand{\Afrak}{\mathfrak{A}}

\newcommand{\abar}{\mathbf{a}}
\newcommand{\bbar}{\mathbf{b}}
\newcommand{\cbar}{\mathbf{c}}
\newcommand{\dbar}{\mathbf{d}}

\newcommand{\cantor}{\{0,1\}^\Nbb}

\newcommand{\Gammao}{\Gamma^o}
\newcommand{\Gammaaa}{\Gamma_{a_0\ldots a_{t-1}}}
\newcommand{\Gammafaa}{\Gamma^f_{a_0\ldots a_{t-1}}}
\newcommand{\spazio}{\;\:}

\newcommand{\labell}[1]{\label{#1}}

\newcommand{\newword}[1]{\textsl{#1}}
\newcommand{\vect}[3]{#1_#2,\ldots ,#1_#3}

\newcommand{\colvect}[3]{(#1_#2\spazio\cdots\spazio #1_#3)^{tr}}

\newcommand{\abs}[1]{\lvert#1\rvert}
\newcommand{\norm}[1]{\lVert#1\rVert_\infty}
\newcommand{\angles}[1]{\langle #1 \rangle}

\DeclareMathOperator{\mmod}{mod}

\theoremstyle{plain}
\newtheorem{theorem}{Theorem}[section]
\newtheorem{lemma}[theorem]{Lemma}
\newtheorem{proposition}[theorem]{Proposition}
\newtheorem{corollary}[theorem]{Corollary}

\theoremstyle{definition}
\newtheorem{definition}[theorem]{Definition}

\newtheorem{remark}[theorem]{Remark}

\usepackage{cite}
\begin{document}

\bibliographystyle{plain}

\sloppy

\title[Kakutani-von Neumann maps]{Kakutani-von Neumann maps on simplexes}

\author[]{Giovanni Panti}
\address{Department of Mathematics\\
University of Udine\\
via delle Scienze 208\\
33100 Udine, Italy}
\email{panti@dimi.uniud.it}

\begin{abstract}
A Kakutani-von Neumann map is the push-forward of the group rotation $(\Zbb_2,+1)$ to a unit simplex $\Gamma$ via an appropriate topological quotient. The usual quotient towards the unit interval $\Gamma=\oi$ is given by the base~$2$ expansion of real numbers, which in turn is induced by the doubling map $x\mapsto 2x\pmod1$. In this paper we replace the doubling map with an $n$-dimensional generalization of the tent map; this allows us to define Kakutani-von Neumann maps in simplexes of arbitrary dimensions. The resulting maps~$K$ are piecewise-linear bijections of $\Gamma$ (not just $\mmod 0$ bijections), whose orbits are all uniformly distributed.
In particular, the dynamical systems $(\Gamma,K)$ are uniquely ergodic, minimal, and have zero topological entropy.
The forward orbit of a certain vertex gives a uniformly distributed enumeration of all points having dyadic coordinates, and this enumeration can be translated via the $n$-dimensional Minkowski function to an enumeration of all rational points in $\Gamma$. In the course of establishing the above results, we introduce a family of $\pm 1$-valued functions on $\Gamma$, constituting an $n$-dimensional analogue of the classical Walsh functions.
\end{abstract}


\thanks{\emph{2000 Math.~Subj.~Class.}: 11K31; 37A45}

\maketitle

\section{Introduction}
We are going to construct certain Kakutani-von Neumann type transformations ---namely, maps which are the push-forward of the adding machine--- on unit simplexes.
Let $\cantor$ be the \newword{Cantor space}, i.e., the space of all countable $0$--$1$ sequences $\abar=a_0a_1a_2\ldots$ endowed with the product topology. It is a key topological fact that every compact metric space ---and in particular every $n$-dimensional simplex--- is a topological quotient (i.e., a continuous surjection) of $\cantor$~\cite[Theorem~3.28]{hockingyou}.

In addition to the topological structure, the Cantor space carries the algebraic structure of the group $\Zbb_2$ of $2$-adic integers, the addition of two group elements $\abar=a_0a_1a_2\ldots$ and $\bbar=b_0b_1b_2\ldots$ being defined componentwise from left to right with carry.
The Haar measure on the compact topological group $\Zbb_2$ is the product measure determined by giving equal mass to the points in $\ooii$.
The \newword{adding machine} $(\Zbb_2,+1)$ is the topological dynamical system induced by addition by $1=10^\infty$ in $\Zbb_2$. It is the prototypical example of a minimal uniquely ergodic dynamical system with zero topological entropy~\cite[IV.1]{CornfeldFomSi82}, \cite[\S15.4]{katokhas95}.
By a \newword{Kakutani-von Neumann transformation} on a unit simplex $\Gamma$ we intend a Borel map from $\Gamma$ to itself which is the push-forward ---in some reasonable sense--- of the adding machine by a topological quotient.

In dimension $1$, the simplest choice for such a quotient is provided by the binary expansion $\varphi(\abar)=\sum a_i2^{-(i+1)}$.
This choice yields the classical Kakutani-von Neumann transformation (also called van der Corput map~\cite[\S5.2.3]{fogg02}) $N:\oi\to\oi$, whose graph is shown in Figure~1.

\begin{figure}[H]
\caption{graph of $N$. The small dots clarify the definition of the function at break points.}
\begin{center}
\includegraphics[height=4.5cm,width=4.5cm]{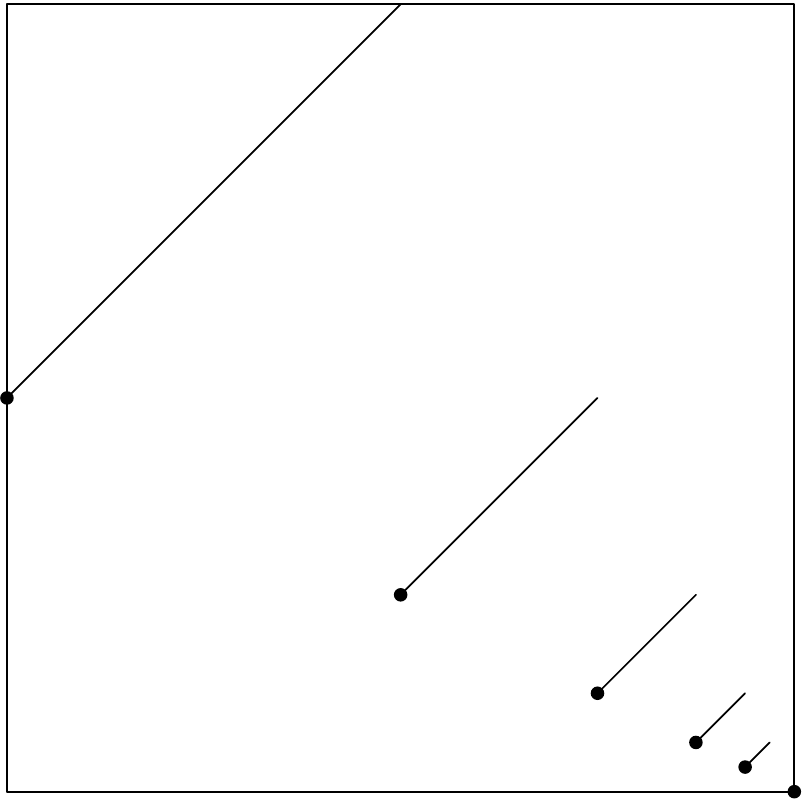}
\end{center}
\end{figure}

The quotient $\varphi$ is induced by the \newword{doubling map} $Dx=2x\pmod2$ on $\oi$. Indeed, $\varphi(\abar)=p$ iff $\abar$ is a \newword{symbolic orbit} for $p$ under $D$ (i.e., for every $t\ge0$, $a_t=0$ implies $D^tp\le1/2$ and $a_t=1$ implies $1/2\le D^tp$).
The point $p$ is \newword{dyadic} iff it belongs to the ring $\Zbb[1/2]$
iff its binary expansion is not unique iff either $p=0$, or $p=1$, or
$D^tp=1/2$ for some $t\ge0$. For points which are not dyadic $N$ is safely defined by $Np=\varphi(\varphi^{-1}(p)+1)$. Things get problematic with dyadic points; namely, the points whose twin expansions are $1^t10^\infty$ and $1^t01^\infty$ (i.e., the points $1/2$, $3/4$, $7/8$, $\ldots$) force us to make choices, expressed by a dot in graphs such as the one in Figure~1. Usually these choices pass unnoticed, since they involve sets of Lebesgue measure $0$. If, however, we are interested in questions such as unique ergodicity, then we cannot discard sets so easily, as they could support invariant measures. 

The dot placing displayed in Figure~1 is the usual one; it corresponds to choosing the finite expansion (i.e., the one ending with $0^\infty$) for each problematic dyadic point. We thus obtain an injective minimal uniquely ergodic map which is not invertible (of course, it is invertible modulo nullsets). 
Note that it is not possible to arrange things so as to obtain an invertible uniquely ergodic map. Indeed, a few moments of reflection on the graph of Figure~1 show that the only possibility of achieving a true bijection is to ``switch dots'' at domain points of the form $(2^i-1)/2^i$, for $i$ odd. In other words, one is forced to map $1/2$ to $1$ (rather than to $1/4$), $7/8$ to $1/4$, and so on. The resulting bijection has finite orbits, such as $0\mapsto 1/2\mapsto 1\mapsto 0$, and hence is not minimal nor uniquely ergodic.
Note also that restricting the domain to the half-open interval $[0,1)$, or gluing $0$ and $1$ together, is of no help, since then~$0$ does not have a counterimage.

The key idea of this paper is to substitute the doubling map $D$ with an $n$-dimensional version of the \newword{tent map} $T$. In dimension~$1$, $T$ is the usual tent map of ergodic theory~\cite[p.~78]{katokhas95}, defined by $T\alpha=2\alpha$ on $[0,1/2]$, and $T\alpha=-2\alpha+2$ on $[1/2,1]$.
In dimension $n\ge1$, $T$ is the transformation introduced in~\cite{panti08} as a linearized version of the \monk\ map. It is a continuous piecewise-linear map on the $n$-dimensional simplex $\Gamma=\{(\vect\alpha1n)\in\Rbb^n:
0\le\alpha_n\le\alpha_{n-1}\le\cdots\le\alpha_1\le1\}$, defined by
$$
T(\alpha_1,\vect \alpha2n)=\begin{cases}
(\alpha_1+\alpha_n,\alpha_1-\alpha_n,\ldots,\alpha_{n-1}-\alpha_n), & \text{if $\alpha_1+\alpha_n\le1$;} \\
(2-\alpha_1-\alpha_n,\alpha_1-\alpha_n,\ldots,\alpha_{n-1}-\alpha_n), & \text{if $\alpha_1+\alpha_n\ge1$.}
\end{cases}
$$
Replacing $D$ with $T$ allows us to construct Kakutani-von Neumann transformations $K$ on unit simplexes of arbitrary dimension. The maps $K$ are piecewise-linear bijections of~$\Gamma$, invertible (not just invertible modulo nullsets), minimal and uniquely ergodic, the unique $K$-invariant Borel measure being the Lebesgue measure.
Actually, we will prove a stronger result, namely that the $K$-orbit of every point in $\Gamma$ is uniformly distributed. We will also prove that the set of points in $\Gamma$ having dyadic coordinates coincides with the forward $K$-orbit of the vertex $(0,0,\ldots,0)$ of $\Gamma$, thus obtaining a uniformly distributed enumeration of these points. 
By conjugating $K$ by the Minkowski question mark function introduced in~\cite{panti08} this yields an enumeration of all rational points in $\Gamma$.
In the course of establishing the above results, we will introduce a family of $\{+1,-1\}$-valued functions on $\Gamma$, constituting an $n$-dimensional analogue of the classical Walsh functions.

\subsubsection*{Acknowledgments} I would like to thank most sincerely Mirko Degli Esposti and Stefano Isola for first suggesting me the possibility of applying the results in~\cite{panti08} to obtain multidimensional versions 
of the Kakutani-von Neumann transformation.
I am also indebted to them for many discussions and several exchanges of emails that greatly helped me in shaping the present work.

\section{Preliminaries}\labell{ref7}

We refer to~\cite{CornfeldFomSi82}, \cite{Walters82}, \cite{katokhas95}, 
for all unexplained notions in topological dynamics and ergodic theory, and to~\cite{rourkesan} for the few needed facts on simplicial complexes.
For the reader's convenience we repeat here the main definitions in~\cite{panti08}. Fix an integer $n\ge1$, and consider the following $(n+1)\times(n+1)$ matrices:
$$
V=\begin{pmatrix}
0 & 1 & 1 & \cdots & 1 & 1 \\
0 & 1 & 1 & \cdots & 1 & 0 \\
0 & 1 & 1 & \cdots & 0 & 0 \\
\vdots & \vdots & \vdots & & \vdots & \vdots \\
0 & 1 & 0 & \cdots & 0 & 0 \\
1 & 1 & 1 & \cdots & 1 & 1
\end{pmatrix},
$$
$$
A_0=\begin{pmatrix}
1 & 0 & 0 & \cdots & 0 & 1 \\
0 & 0 & 0 & \cdots & 0 & 1 \\
0 & 1 & 0 & \cdots & 0 & 0 \\
\vdots & \vdots & \vdots & & \vdots & \vdots \\
0 & 0 & 0 & \cdots & 0 & 0 \\
0 & 0 & 0 & \cdots & 1 & 0
\end{pmatrix},
\qquad
A_1=\begin{pmatrix}
0 & 0 & 0 & \cdots & 0 & 1 \\
1 & 0 & 0 & \cdots & 0 & 1 \\
0 & 1 & 0 & \cdots & 0 & 0 \\
\vdots & \vdots & \vdots & & \vdots & \vdots \\
0 & 0 & 0 & \cdots & 0 & 0 \\
0 & 0 & 0 & \cdots & 1 & 0
\end{pmatrix}.
$$
More precisely: all entries of $V$ are $0$, except those in position $ij$, with either ($i=n+1$) or ($j\ge2$ and $i+j\le n+2$), that have value $1$. All entries of $A_0$ and $A_1$ are $0$, except $(A_0)_{11}$, $(A_1)_{21}$, and all elements in position $1(n+1)$, $2(n+1)$, $(j+1)j$, for $2\le j\le n$, that have value $1$.

For $a=0,1$, let $B_a$ be the matrix which is identical to $A_a$, except for the last column, where the two entries $1$ are replaced by $1/2$. Let $C_a=VA_aV^{-1}$, and 
$D_a=VB_aV^{-1}$. We write points $(\vect\alpha1n)\in\Rbb^n$ as column vectors in projective coordinates $(\alpha_1\cdots\alpha_n\,1)^{tr}$, and we let $\Gamma$ be the $n$-dimensional simplex in $\Rbb^n$ whose $i$th vertex $v_i$ is given by the $(i+1)$th column of $V$, for $0\le i\le n$;
in affine coordinates,
$\Gamma=\{(\vect\alpha1n)\in\Rbb^n:
0\le\alpha_n\le\alpha_{n-1}\le\cdots\le\alpha_1\le1\}$.

The matrix $B_1$ is column-stochastic and primitive (i.e., some power is strictly positive). By the Perron-Frobenius theory its conjugate $D_1$ has exactly one eigenvector $(\alpha_1\cdots\alpha_n\,1)^{tr}$ such that the point $(\vect{\alpha}1n)$, which we denote by~$v_{-1}$, is in~$\Gamma$. The corresponding eigenvalue is $1$, and hence $\vect\alpha 1n$ are rational numbers; for example, for $n=1,2,3$ we have $v_{-1}=2/3,(4/5,2/5),(6/7,4/7,2/7)$, respectively.
For $a=0,1$, let $\tau_a:\Gamma\to\Gamma$ be the affine map determined by $D_a$, namely $\tau_a(\vect\alpha1n)=(\vect\beta1n)$ iff $D_a(\alpha_1\cdots\alpha_n\,1)^{tr}=(\beta_1\cdots\beta_n\,1)^{tr}$.
As proved in~\cite{panti08}, $\tau_0$ and $\tau_1$ are the two inverse branches of the tent map $T$ defined in the Introduction.
Let $\Gammao$ be the set of all points $\sum x_iv_i\in\Gamma$ such that $\vect x0n\ge0$, $\sum x_i=1$, and $x_0>0$; one easily computes that $\Gamma\setminus\Gammao=\{(\vect\alpha 1n)\in\Gamma:\alpha_1=1\}$.

\begin{proposition}
The\labell{ref5} sets $\{v_{-1}\}$, $\tau_0\Gammao$, $\tau_1\tau_0\Gammao$, $\tau_1^2\tau_0\Gammao$, $\tau_1^3\tau_0\Gammao$, $\ldots$ constitute a partition of $\Gamma$; the same holds for the sets $\{v_0\}$, $\tau_1\Gammao$, $\tau_0\tau_1\Gammao$, $\tau_0^2\tau_1\Gammao$, $\tau_0^3\tau_1\Gammao,\ldots$.
\end{proposition}
\begin{proof}
We claim that for every $k\ge0$ the sets $\tau_0\Gammao, \tau_1\tau_0\Gammao, \ldots, \tau_1^k\tau_0\Gammao, \tau_1^{k+1}\Gamma$ constitute a partition of $\Gamma$. This is true for $k=0$ since, as proved in~\cite{panti08}, $\Gamma=\tau_0\Gamma\cup \tau_1\Gamma$ and $\tau_0\Gamma\setminus \tau_1\Gamma=\tau_0\Gammao$. Assume the statement true for $k$. Since $\tau_1:\Gamma\to \tau_1\Gamma$ is a bijection, $\tau_1^{k+1}:\Gamma\to \tau_1^{k+1}\Gamma$ is a bijection as well, and hence the partition $\tau_0\Gammao,\tau_1\Gamma$ of $\Gamma$ induces a partition $\tau_1^{k+1}\tau_0\Gammao,\tau_1^{k+2}\Gamma$ of $\tau_1^{k+1}\Gamma$; this settles our claim. The first statement now follows readily since, by construction, $\bigcap_{k\ge1}\tau_1^k\Gamma=\{v_{-1}\}$. The proof of our second statement is analogous, using $\tau_1\Gammao,\tau_0\Gamma$ as base partition and observing that $\bigcap_{k\ge1}\tau_0^k\Gamma=\{v_0\}$.
\end{proof}

We can now define our Kakutani-von Neumann transformation $K:\Gamma\to\Gamma$. First, we set $Kv_{-1}=v_0$; second, for every $k\ge0$ we have bijections
$$
\begin{xy}
\xymatrix{
\tau_1^k\tau_0\Gammao & \Gammao \ar[l]_-{\tau_1^k\tau_0} \ar[r]^-{\tau_0^k\tau_1} & \tau_0^k\tau_1\Gammao,
}
\end{xy}
$$
and we define $K=\tau_0^k\tau_1(\tau_1^k\tau_0)^{-1}$ on $\tau_1^k\tau_0\Gammao$.

A couple of pictures may be helpful. In Figure~2 we draw the graph of $K$ for $n=1$ and $\Gamma=\oi$. In Figure~3 we draw the partitions $\{\tau_0\Gammao$,  $\tau_1\tau_0\Gammao$, $\tau_1^2\tau_0\Gammao$, $\tau_1^3\tau_0\Gammao$, $\tau_1^4\tau_0\Gammao$, $\tau_1^5\Gamma\}$ and
$\{\tau_1\Gammao$, $\tau_0\tau_1\Gammao$, $\tau_0^2\tau_1\Gammao$, $\tau_0^3\tau_1\Gammao$, $\tau_0^4\tau_1\Gammao$, $\tau_0^5\Gamma$\}, for $n=2$; here $\Gamma=\{(\alpha,\beta):0\le\beta\le\alpha\le1\}$.
\begin{figure}[H]
\caption{graph of $K$ in dimension $1$.}
\begin{center}
\includegraphics[height=4.5cm,width=4.5cm]{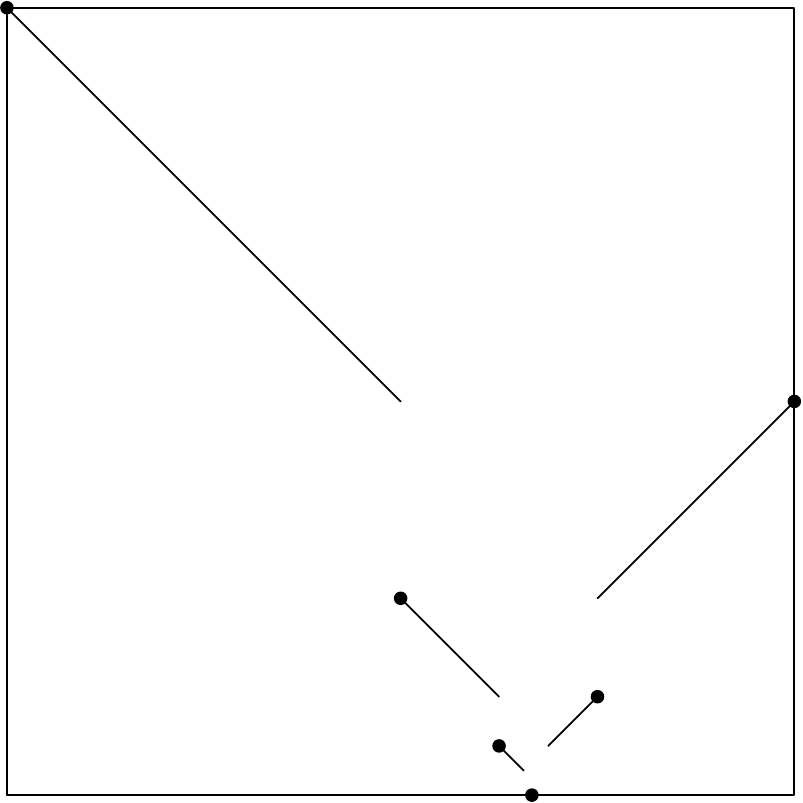}
\end{center}
\end{figure}
\begin{figure}[H]
\caption{partitions of $\Gamma$, in dimension~$2$.}
\begin{center}
\includegraphics[height=4.5cm,width=4.5cm]{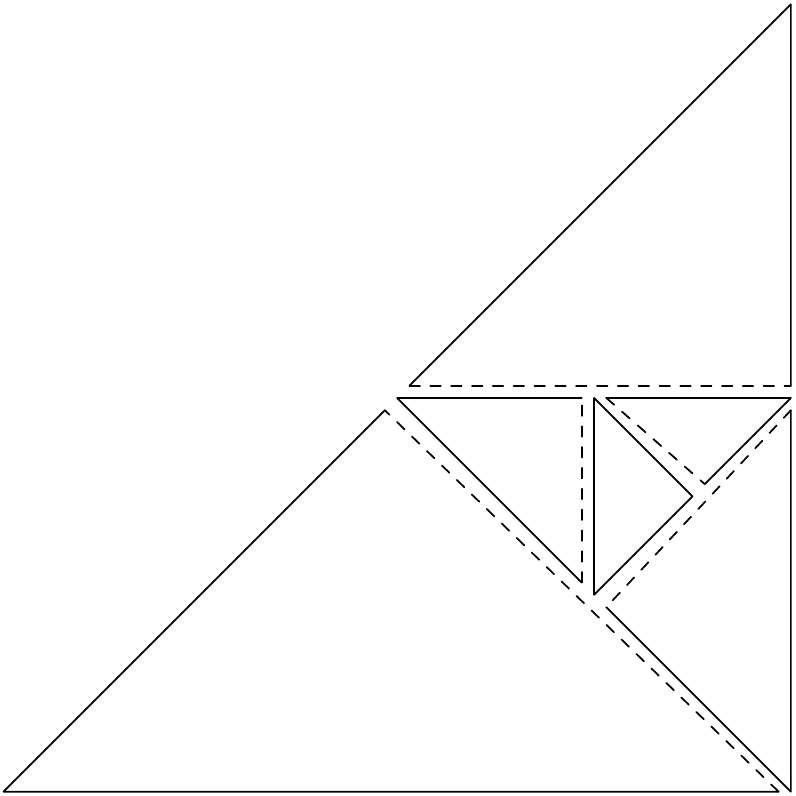}
\qquad
\includegraphics[height=4.5cm,width=4.5cm]{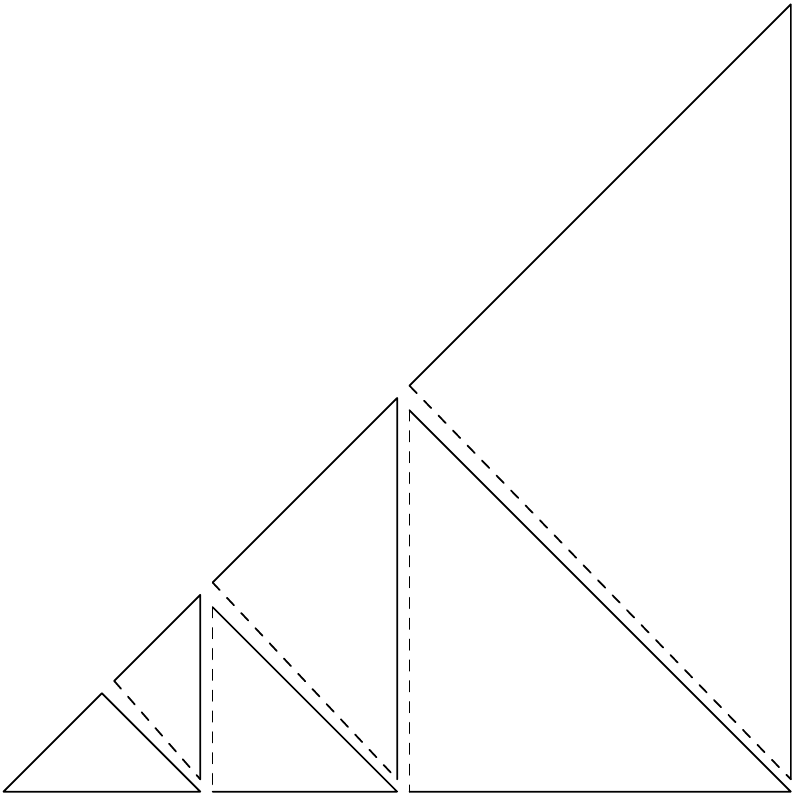}
\end{center}
\end{figure}

A good way of comparing the classical Kakutani-von Neumann map $N:\oi\to\oi$ with the 
$1$-dimensional version of our $K$ is by noting that $N$ is definable by a cut-and-stack procedure~\cite{friedman92}. At stage $1$, the interval $\oi$ is cut into two equal pieces, and the right-hand piece $[1/2,1]$ is stacked on top of the left-hand one $[0,1/2)$ in an orientation-preserving way. At stage $k+1$, the $2^k$ layers of stage $k$ are all cut in two equal pieces, and the resulting right-hand pieces are stacked on top of the left-hand ones, again in an orientation-preserving way. At stage $k$ a map $N_k$ is defined as the map that moves every point of every layer (except the top one) to the point immediately above in the next layer. A clear limiting process defines then $N$ from the partial maps $N_k$, and shows that $N$ is a bijection modulo $\lambda$-nullsets (we always use $\lambda$ to denote the Lebesgue measure, normalized so that $\lambda(\Gamma)=1$).
Now, if we are happy to neglect $\lambda$-nullsets (namely, the endpoints of the cut intervals), then the same cut-and-stack procedure defines our $K$, except that now at each stage the right-hand half-intervals must be put on top of the left-hand ones \emph{in an orientation-reversing way}. We leave to the reader the straightforward verification that this construction is correct, i.e., defines a map $\lambda$-everywhere identical to $K$.
We stress however that $K$ is a true bijection on~$\Gamma$ in every dimension $n\ge1$, not just a mod $0$ one.

\section{Coding points}\labell{ref11}

For $t\ge1$ and $\vect a0{{t-1}}\in\ooii$, let $\Gammaaa$ be the simplex
$\tau_{a_0}\tau_{a_1}\cdots\tau_{a_{t-1}}\Gamma$. We have the identity
\begin{equation}\label{eqref1}
\Gammaaa=
\Gamma_{a_0}\cap T^{-1}\Gamma_{a_1}\cap T^{-2}\Gamma_{a_2}
\cap\cdots\cap T^{-(t-1)}\Gamma_{a_{t-1}}.
\end{equation}
Moreover, for a fixed $t$ the set of all faces of all $n$-dimensional simplexes $\Gammaaa$ forms a simplicial complex $\Bcal_t$ whose support is $\Gamma$. Each complex in the chain $\Bcal_0, \Bcal_1, \Bcal_2, \ldots$ refines the preceding ones; see~\cite[Proposition~2.2]{panti08} for the above statements.

By~\cite[Lemma~2.3]{panti08}, for every point $\abar=a_0a_1a_2\ldots$ of $\Zbb_2$ the set 
$\bigcap\{\Gammaaa:t\ge0\}$ is a singleton $\{p\}$, and the map $\upsilon:\Zbb_2\to\Gamma$ defined by $\upsilon(\abar)=p$ is a topological quotient.
The $\upsilon$-counterimage of $\Gammaaa$ is the cylinder $[\vect a0{{t-1}}]$ whose elements are all sequences $\bbar$ such that $b_i=a_i$ for every $0\le i<t$.

\begin{definition}
If $\upsilon(\abar)=p$, we say that $\abar$ is a \newword{symbolic orbit} of $p$ under $T$. 
By~\eqref{eqref1}, this amounts to saying that $T^tp\in\Gamma_{a_t}$ for every $t\ge0$.
Since $\Gamma_0\cap\Gamma_1\not=\emptyset$, the symbolic orbit is not unique. The simplex $\Gammaaa$ is the set of points that have a symbolic orbit beginning with $a_0\ldots a_{t-1}$.
Let now $\cbar\in\{0,1,*\}^\Nbb$; we say that $\cbar$ is the 
$*$-orbit of $p$ if $c_t$ equals $0$, $1$, or $*$, whenever $T^tp$ belongs to $\Gamma_0\setminus\Gamma_1$, $\Gamma_1\setminus\Gamma_0$, or $\Gamma_0\cap\Gamma_1$, respectively. Every point has a unique $*$-orbit, and replacing the symbols~$*$ in it by $0$ or $1$, arbitrarily, we obtain its various symbolic orbits.
The \newword{final orbit} of $p$ is its unique symbolic orbit $\abar$ such that $T^tp\in\Gamma_0\cap\Gamma_1$ implies $a_t=1$; it is obtained from the $*$-orbit by replacing all symbols $*$ by $1$.
\end{definition}

\begin{theorem}
Every\labell{ref1} $p\in\Gamma$ has at most $2^{n(n+1)/2}$ symbolic orbits.
\end{theorem}
\begin{proof}
We will show that the $*$-orbit of $p$ contains at most $1+2+3+\cdots +n$ symbols~$*$. Consider the following cones in $\Rbb^{n+1}$:
\begin{align*}
C&=\{\colvect\beta1{{n+1}}:\beta_i\ge0\text{ for every $i$, and }
\beta_i>0\text{ for at least one $i$}\};\\
C_0&=\{\colvect\beta1{{n+1}}\in C:\beta_1\ge\beta_2\};\\
C_1&=\{\colvect\beta1{{n+1}}\in C:\beta_1\le\beta_2\}.
\end{align*}
Let $\Mcal:C\to C$ be defined by
$$
\Mcal\colvect\beta1{{n+1}}=
\begin{cases}
(\beta_1-\beta_2\spazio\beta_3\spazio\beta_4\spazio\cdots\spazio\beta_{n+1}\spazio\beta_2)^{tr},
&\text{if $\beta_1\ge\beta_2$};\\
(\beta_2-\beta_1\spazio\beta_3\spazio\beta_4\spazio\cdots\spazio\beta_{n+1}\spazio\beta_1)^{tr},
&\text{if $\beta_1\le\beta_2$}.
\end{cases}
$$
By explicit computation one sees that $\Mcal$ is induced by left multiplication by the inverses of the matrices $A_0, A_1$ defined in Section~\ref{ref7}.
By~\cite[p.~250]{panti08}, $\Mcal$ is the projective version of a conjugate of the \monk\ map. Let $P=V^{-1}\Phi^{-1}(p)$; here $\Phi$ is the Minkowski question mark function in~\cite[Theorem~2.1]{panti08}, and $\Phi^{-1}(p)$ is written as a column vector in projective coordinates.
By~\cite[Lemma~2.5]{panti08} we have $T^tp\in\Gamma_a$ iff $\Mcal^tP\in C_a$,
for $a=0,1$. Hence the set of $T$-symbolic orbits of $p$ coincides with the set of $\Mcal$-symbolic orbits (taken w.r.t.~the partition $C_0,C_1$) of $P$. 

For simplicity's sake, write $\Mcal^t P=P^t=\colvect{\beta^t}1{{n+1}}$, and refer to a time $t$ at which $P^t\in C_0\cap C_1$ as a \newword{hitting time}. The hitting time $t$ is \newword{principal} if $\beta_1=\beta_2\not=0$, and is \newword{secondary} if $\beta_1=\beta_2=0$. Let $0\le z(t)\le n$ be the number of symbols~$0$ appearing in $\colvect{\beta^t}1{{n+1}}$. By the definition of $\Mcal$ we always have $z(t+1)=z(t)$, except when $t$ is a primary hitting time, in which case $z(t+1)=z(t)+1$. Therefore, following the $\Mcal$-orbit of $P$ we will encounter at most $n$ primary hitting times. Say that $r\ge0$ is one of these, and let $P^{r+1}=(0\spazio\cdots\spazio 0\spazio\beta^{r+1}_{m+1}\spazio\cdots\spazio\beta^{r+1}_{n+1})^{tr}$, with $\beta^{r+1}_{m+1}\not=0$ and $m\ge1$.
Then $r+1,r+2,\ldots,r+m-1$ are all secondary hitting times (there are none such if $m=1$), and
$P^{r+m}=(0\spazio\beta^{r+1}_{m+1}\spazio\cdots\spazio
\beta^{r+1}_{n+1}\spazio 0\spazio\cdots\spazio 0)^{tr}$.
Since $\beta^{r+1}_{m+1}\not=0$, a moment of reflection shows that the first hitting time $t> r+m$ must be a primary one. Since $m\le z(r+1)=z(r)+1$, we have $m-1\le z(r)$; in other words, every primary hitting time $r$ is followed by at most $z(r)$ secondary hitting times.
By the same argument, if the first hitting time $r$ is secondary, then it is followed by at most $z(r)-2$ secondary hitting times.
Since the number of symbols~$0$ in the coordinates of points in $C$ varies from $0$ to $n$, it is clear that the number of hitting times along the $\Mcal$-orbit of $P$ is bounded by the triangular number $1+2+3+\cdots+n=n(n+1)/2$.
\end{proof}

\begin{remark}
The\labell{ref4} set of admissible $*$-orbits depends on the dimension $n$. For example, consider the sequence $\cbar=0010{*}{*}10^\infty$. By Theorem~\ref{ref1}, $\cbar$ is not realizable as a $*$-orbit in dimension~$1$. Let $n\ge2$; the only point in $\Gamma$ whose $*$-orbit may possibly equal $\cbar$ is $p=\tau_0^2\tau_1\tau_o\tau_1^3v_0$ (see Lemma~\ref{ref6} below). By actual computation, one checks that in dimension~$2$ $p$ equals $(3/8,3/8)$, whose $*$-orbit is indeed $\cbar$. On the other hand, in dimension~$3$ $p$ equals $(3/4,1/4,1/4)$ and has $*$-orbit ${*}{*}1{*}{*}{*}10^\infty$, while in dimension~$4$ $p=(1/2,1/2,0,0)$, whose $*$-orbit is $00{*}{*}{*}{*}10^\infty$.
\end{remark}

By saying that $\abar$ is a \newword{final orbit} we mean that $\abar$ is the final orbit of some point, necessarily the point $\upsilon(\abar)$.

\begin{lemma}
The\labell{ref6} following statements hold.
\begin{itemize}
\item[(i)] If $\bbar$ is a symbolic orbit of $p$, then $a\bbar$ is a symbolic orbit of $\tau_ap$, for $a=0,1$.
\item[(ii)]
\begin{gather*}
\tau_0\Gammao=\Gamma\setminus\Gamma_1=\Gamma_0\setminus\Gamma_1=\{p\in\Gamma:\text{ the final orbit of $p$ begins with $0$}\},\\
\tau_1\Gamma=\Gamma_1=\{p\in\Gamma:\text{ the final orbit of $p$ begins with $1$}\}.
\end{gather*}
\item[(iii)] If $p\in\Gammao$ has final orbit $\bbar$, then $\tau_ap$ has final orbit $a\bbar$, for $a=0,1$.
\item[(iv)] If $p\in\Gamma\setminus\Gammao$ has final orbit $\bbar$, then $\tau_0p=\tau_1p\in\Gamma_0\cap\Gamma_1$ and has final orbit $1\bbar$.
\item[(v)] If $w\bbar$ is a final orbit, with $w$ a finite $\ooii$-word, then $1^t\bbar$ is a final orbit for every $t\ge0$.
\item[(vi)] $0^\infty$ and $1^\infty$ are the unique symbolic orbits of the points $v_0$ and $v_{-1}$, respectively.
\end{itemize}
\end{lemma}
\begin{proof}
Since $T\tau_a=id_\Gamma$, (i) is clear.
The second identity in (ii) is a special case of~\eqref{eqref1}, and the first follows by taking complements and looking at the proof of Proposition~\ref{ref5}. (iii) follows from (ii), and (iv) from (i) and the observation that $\tau_0=\tau_1$ on $\Gamma\setminus\Gammao$.
(v) If $w\bbar$ is the final orbit of $p$, then $1^t\bbar$ is the final orbit of $\tau_1^tT^{\abs w}p$, by (iii) and (iv).
We obtain (vi) by noting that $v_0$ and $v_{-1}$ are the only fixed points for $\tau_0$ and $\tau_1$, respectively.
\end{proof}

In the following main Theorem~\ref{ref2} we will give an alternate description of $K$. As corollaries, we will obtain that the forward $K$-orbit of $v_0$ contains precisely all dyadic points in~$\Gamma$ (a point is \newword{dyadic} if all its coordinates are dyadic), each such point appearing exactly once. We will also obtain that the measure-preserving system $(\Gamma, \lambda, K)$ is metrically isomorphic to the adding machine $(\Zbb_2,\text{Haar measure},+1)$.

Let $\abar\in\Zbb_2$ and, for each $t\in\Zbb$, let $\abar^t=\abar+t$. Consider the doubly infinite bilateral orbit $\Afrak$ of $\abar$ under $+1$, namely $\Afrak=\ldots,\abar^{-2},\abar^{-1},\abar^0=\abar,\abar^1,\abar^2,\ldots$. If $\abar$ ends with either $0^\infty$ or $1^\infty$, then the elements of $\Afrak$ are all sequences $\bbar$ such that $\bbar$ ends with either $0^\infty$ or $1^\infty$. Otherwise, the elements of $\Afrak$ are all sequences $\bbar$ that \newword{have the same tail} of $\abar$ (i.e., such that there exists $t\ge0$, depending on $\bbar$, with $a_r=b_r$ for every $r\ge t$). Let $p\in\Gamma$ be such that $\Afrak$ contains at least one symbolic orbit of~$p$. From Theorem~\ref{ref1} and the above characterization of the elements of $\Afrak$, this implies that $\Afrak$ contains all symbolic orbits of $p$. For every such $p$, remove from $\Afrak$ all these symbolic orbits \emph{except the last one} (this last one is precisely the final orbit of $p$, whence the name). 
Write $\Afrak'=\ldots,\abar^{r_{-2}},\abar^{r_{-1}},\abar^{r_0},\abar^{r_1},\abar^{r_2},\ldots$ for the resulting pruned sequence. By Theorem~\ref{ref1}, $\Afrak'$ is surely infinite to the right, and the following Theorem~\ref{ref2} implies, since $K$ is a bijection, that it is infinite to the left as well.

\begin{theorem}
Let\labell{ref2} $\abar^{r_k}\in\Afrak'$ be the final orbit of $p$. Then the final orbit of $Kp$ is $\abar^{r_{k+1}}$.
\end{theorem}
\begin{proof}
Without loss of generality $\abar^0$ is a final orbit and $r_k=r_0=0$. If $\abar^0=1^\infty$, then $\abar^0$ is the final orbit of $v_{-1}$ and $\abar^1=0^\infty$ is the final orbit of $v_0=Kv_{-1}$. By the definition of $\Afrak'$ we have $\abar^{r_1}=\abar^1$ and we are through. Otherwise, write uniquely $\abar^0=1^t0\bbar$ for a certain $t\ge0$, and let $q=T^{t+1}p$; then $q$ has final orbit $\bbar$. We claim that $q\in\Gammao$. Indeed, the point $T^tp$ has final orbit $0\bbar$, and hence is in $\tau_0\Gammao$ by Lemma~\ref{ref6}(ii). Since $T\tau_0$ is the identity map, we have $q=TT^tp\in T\tau_0\Gammao=\Gammao$. By (iii) and (iv) of the same Lemma, the point $\tau_1^t\tau_0q$ has final orbit $\abar^0$, and hence is $p$. By the definition of $K$, we have $Kp=\tau_0^t\tau_1q$, and by (i) one of the symbolic orbits of $Kp$ (not necessarily the final one) is $0^t1\bbar=\abar^1$. We will show that $\upsilon(\abar^1)=\upsilon(\abar^{r_1})$, thus concluding the proof.

Define numbers $t_1,t_2,\ldots$ in the following inductive way:
\begin{itemize}
\item $t_1$ is the maximum number $\ge0$ and $\le t$ such that $0^{t_1}1\bbar$ is a final orbit.
\item Assume $t_i$ has been defined. If $s_i=(t+1)-\sum_{j=1}^i(t_j+1)=0$, then the procedure stops at $t_i$. Otherwise $t_{i+1}$ is the maximum number $\ge0$ and $\le t-s_i$ such that $0^{t_{i+1}}10^{t_i}1\cdots 0^{t_2}10^{t_1}1\bbar$ is a final orbit.
\end{itemize}
After finitely many steps the procedure stops, say at $t_d$. For $1\le i\le d$, let now $u_i$ be the integer whose binary expansion (written from left to right) is
$0^{s_i}0^{t_i}10^{t_{i-1}}1\cdots 0^{t_1}1$; this expansion contains exactly $t+1$ digits. Let $l_i=u_i-(2^t-1)$; observe that $2^t-1$ has binary expansion $1^t0$. We then have $1=l_1<l_2<\cdots<l_d$. By construction,
$\abar^{l_i}=0^{s_i}0^{t_i}10^{t_{i-1}}1\cdots 0^{t_1}1\bbar$, which is a final orbit iff $i=d$ (note that $s_d=0$). Now, let $1\le i\le d$ be maximum such that $\upsilon(\abar^{l_1})=\upsilon(\abar^{l_2})=\cdots=\upsilon(\abar^{l_i})$ and for every $l_1\le l<l_i$ the orbit $\abar_l$ is not final. We shall show that $i=d$, so that $l_d=r_1$ (since $\abar^{l_d}$ is a final orbit) and $\upsilon(\abar^1)=\upsilon(\abar^{r_1})$, thus concluding the proof.

Assume by contradiction $i<d$. Then $0^{s_i}0^{t_i}10^{t_{i-1}}1\cdots 0^{t_1}1\bbar$ is not final. By the definition of $t_i$, the orbit
$0^{t_i}10^{t_{i-1}}1\cdots 0^{t_1}1\bbar$ is final, while
$00^{t_i}10^{t_{i-1}}1\cdots 0^{t_1}1\bbar$ is not.
Therefore, by Lemma~\ref{ref6}(iii), the point
$z=\upsilon(0^{t_i}10^{t_{i-1}}1\cdots 0^{t_1}1\bbar)$ does not belong to $\Gammao$. Applying (iv), we obtain
$\upsilon(\abar^{l_i})=\tau_0^{s_i}z=\tau_0^{s_i-1}\tau_1z=
\upsilon(\abar^{l_{i+1}})$.
Moreover, for every $l_i\le l<l_{i+1}$ the orbit $\abar^l$ is not final, since it has $00^{t_i}10^{t_{i-1}}1\cdots 0^{t_1}1\bbar$ as a tail. This contradicts the maximality of $i$ and concludes the proof.
\end{proof}

\begin{corollary}
Two points belong to the same $K$-orbit iff either their final orbits have the same tail, or one of them has tail $0^\infty$ and the other $1^\infty$ (one can use equivalently any symbolic orbit, or the $*$-orbit).
\end{corollary}
\begin{proof}
Immediate from Theorem~\ref{ref2}.
\end{proof}

\begin{corollary}
The\labell{ref18} forward $K$-orbit of $v_0$ constitutes an enumeration without repetitions of all dyadic points in $\Gamma$.
\end{corollary}
\begin{proof}
By~\cite[Theorem~3.5(ii)]{panti08} the dyadic points in $\Gamma$ are exactly the points whose symbolic orbits have tail $0^\infty$. Since $v_0$ has symbolic orbit $0^\infty$, if we set $\abar^0=0^\infty$ and construct $\Afrak'$ as above, then the points in $\Afrak'$ to the right of $\abar^0$ are exactly the final orbits of the dyadic points.
\end{proof}

Since the proof of Corollary~\ref{ref18} makes crucial use of~\cite[Theorem~3.5]{panti08}, we take this opportunity to correct the annoying ---albeit apparent--- typo in the statement~(iii) of the above reference, where ``iff'' must be read ``if''.

\begin{corollary}
The\labell{ref17} map $K$ preserves the Lebesgue measure $\lambda$. The 
measure-preserving system $(\Gamma,\lambda,K)$ is metrically isomorphic to the adding machine $(\Zbb_2,\mu,+1)$, where $\mu$ is the Haar measure on $\Zbb_2$.
\end{corollary}
\begin{proof}
The matrix $E_k=D_0^kD_1D_0^{-1}D_1^{-k}$ induces $K$ on $\tau_1^k\tau_0\Gammao$.
We have $\abs{det(D_0)}=\abs{det(D_1)}=1/2$, and hence $\abs{det(E_k)}=1$.
The row vector $(1\cdots1\,1)$ is a left eigenvector for the column-stochastic matrices $B_0$ and $B_1$, and therefore $(1\cdots1\,1)V^{-1}=(0\cdots0\,1)$ is a left eigenvector for $D_0$ and $D_1$, with corresponding eigenvalue $1$. This implies that the last row of $E_k$ is $(0\cdots0\,1)$, and hence the Jacobian matrix of $K$ on $\tau_1^k\tau_0\Gammao$ is the $n\times n$ minor given by the first $n$ rows and columns of $E_k$. Such a matrix has determinant of absolute value $1$, and hence $K$ preserves the Lebesgue measure.

For every $t>0$ and every $(\vect a0{{t-1}})\in\ooii^t$ the Lebesgue measure 
$2^{-t}$ of $\Gammaaa$ coincides with the Haar measure of its $\upsilon$-counterimage $[\vect a0{{t-1}}]$.
Since the sets 
$\Gammaaa$ generate the Borel $\sigma$-algebra of $\Gamma$, the push-forward of $\mu$ via $\upsilon$ is $\lambda$.
Let $Y=\{p\in\Gamma:$ the symbolic orbit of $p$ is not unique$\}$. Then $Y$ is a $\lambda$-nullset, since it is contained in the $\lambda$-nullset
$\bigcup_{t\ge0}T^{-t}[\Gamma_0\cap\Gamma_1]$. Therefore $X=\Gamma\setminus
\bigcup_{t\in\Zbb}K^tY$ has full Lebesgue measure, and $\upsilon^{-1}X$ has full Haar measure. By construction, $\upsilon$ is a bijection between $\upsilon^{-1}X$ and $X$. If $\abar\in\upsilon^{-1}X$ and $p=\upsilon(\abar)$, then all the elements in the $K$-orbit of $p$ have a unique symbolic orbit. Hence the sequence $\Afrak$ constructed from $\abar$ coincides with $\Afrak'$, and by Theorem~\ref{ref2} $\upsilon(\abar+1)=K\upsilon(\abar)$.
\end{proof}

\section{$T$-Walsh functions on simplexes}\labell{ref15}

In this section we will prove that every $K$-orbit is $\lambda$-uniformly distributed.
Remember~\cite[Definition~III.1.1]{kuipersnie74} that a countable sequence $\{p^i\}_{i\in\Nbb}$ in a compact metric space $X$ endowed with a Borel probability measure $\nu$ is
\newword{$\nu$-uniformly distributed} if for every continuous function $f:X\to\Rbb$ we have
\begin{equation}\label{eqref2}
\lim_{k\to\infty}\frac{1}{k}\sum_{i=0}^{k-1}f(p^i)=
\int_X f\,d\nu.
\end{equation}

Let $\iota:\Gamma\to\cantor$ be the function that associates to a point its final orbit. It is a right inverse to $\upsilon$, so it is injective. For $t\ge1$ and 
$\vect a0{{t-1}}\in\ooii$, let 
$\Gammafaa=\iota^{-1}[\vect a0{{t-1}}]=\{p\in\Gamma:$ the final orbit of $p$ begins with 
$a_0\ldots a_{t-1}\}$. For a fixed $t$, the collection of all sets $\Gammafaa$ is a true partition of~$\Gamma$, not just a $\mmod 0$ one.

\begin{lemma}
The set $\Gammafaa$ is\labell{ref8} obtained from the simplex $\Gammaaa$ by removing some of its proper faces. In particular, $\Gammafaa$ and $\Gammaaa$ have the same Lebesgue measure $2^{-t}$.
The embedding $\iota$ is Borel, not continuous, and has dense range.
\end{lemma}
\begin{proof}
Let $w=a_0\ldots a_{t-1}$ be a finite word, let $0\le i<t$, and let $w(i,1)$ be the word obtained from $w$ by setting $a_i=1$. We immediately have from the definitions
$$
\Gamma^f_w=\Gamma_w\setminus\bigl(\bigcup\{\Gamma_{w(i,1)}:a_i=0\}\bigr).
$$
Since $\Bcal_t$ is a simplicial complex, if $a_i=0$ then $\Gamma_w\cap\Gamma_{w(i,1)}$ is a proper face (possibly empty) of $\Gamma_w$, and our first two statements follow. As the clopen cylinders $[\vect a0{{t-1}}]$ generate the topology of the Cantor space, $\iota$ is Borel; it is not continuous since $\iota^{-1}[0]=\Gamma^f_0=\Gamma_0\setminus\Gamma_1$ is not closed. Finally, $\iota$ has dense range because no $\Gammafaa$ is empty.
\end{proof}

We now change the group structure on the Cantor space, by endowing it with the structure of the direct product $Z_2^\Nbb$ of countably many copies of the two-element group $Z_2$ (no blackboard type here). Addition of two elements $\abar$ and $\bbar$ is done componentwise without carry, the topology and the Haar measure being the same as those of $\Zbb_2$.

Just as a Kakutani-von Neumann map is the push-forward of the adding machine $(\Zbb_2,+1)$ by a topological quotient, a \newword{Walsh function} is the push-forward of a character $\chi$ of $Z_2^\Nbb$ by the same topological quotient; in our case, it is simply the composition $\chi\iota$. By the Pontrjagin duality, the character group of $Z_2^\Nbb$ is the direct sum $\bigoplus^\Nbb Z_2$ with the discrete topology; since $Z_2^\Nbb$ has exponent $2$, each character has range in $\{+1,-1\}$. By associating the element $(d_0,d_1,\ldots,d_{t-1},0,0,\ldots)$ (with $d_{t-1}=1$) to $m=\sum_{i=0}^{t-1}d_i2^i$, we identify $\bigoplus^\Nbb Z_2$ with $\Nbb$; under this identification the sum of $m$ and $l$ is the natural number $m\oplus l$ whose $i$th binary digit is the $\mmod2$ sum of the $i$th binary digits of $m$ and $l$. In short, we obtain the following definition.

\begin{definition}
Given\labell{ref9} $m\in\Nbb$ whose binary expansion is $m=\sum_{i=0}^{t-1}d_i2^i$ (with $d_{t-1}=1$) and $p\in\Gammafaa$, write $\angles{m,p}$ for $d_0a_0+d_1a_1+\cdots+d_{t-1}a_{t-1}\pmod2$.
The $m$th $T$-\newword{Walsh function} $u_m=\chi_m\iota:\Gamma\to\{+1,-1\}$ is defined by $u_m(p)=(-1)^{\angles{m,p}}$; we have $u_0=\one_\Gamma$.
For $m\ge1$ the \newword{level} of $u_m$ is $t\ge1$, and the set $U_t$ of $T$-Walsh functions of level $t$ has cardinality $2^{t-1}$.
\end{definition}

The name $T$-Walsh refers to the r\^ole of the map $T$ ---that stays on the stage through the $\iota$ embedding--- in generating the group $U=\{u_0\}\cup\bigcup_{t\ge1}U_t$; see Proposition~\ref{ref10}(i).
Write $r=u_1$ for the first $T$-Walsh function; then $r(p)$ equals $+1$ or $-1$ according whether $p$ belongs to $\Gamma^f_0$ or not. The functions $r,rT,rT^2,\ldots$ correspond to the classical Rademacher functions; see Remark~\ref{ref16}.

\begin{proposition}
Let\labell{ref10} $m=\sum_{i=0}^{t-1}d_i2^i$ be as above; then the following statements hold.
\begin{itemize}
\item[(i)]
$$
u_m=\prod_{i=0}^{t-1}(rT^i)^{d_i};
$$
this expression is unique and includes the case $m=0$, since an empty product equals $1$ by definition.
\item[(ii)] $u_m\cdot u_l=u_{m\oplus l}$, and hence the group algebra $\Rbb[U]$ coincides with the $\Rbb$-span of $U$.
\item[(iii)] The family $U$ forms an orthonormal set in $L_2(\Gamma,\Rbb)$.
\item[(iv)] The closure of $\Rbb[U]$ in $L_\infty(\Gamma,\Rbb)$ contains all continuous real-valued functions on $\Gamma$.
\end{itemize}
\end{proposition}

\begin{proof}
(i) Under the identification of the character group of $Z_2^\Nbb$ with $\bigoplus^\Nbb Z_2$, the function $r=\chi_1\iota$ corresponds to $(1,0,0,\ldots)$. Applying $T$ amounts to shifting the final orbit of a point by one, and therefore $rT^i$ corresponds to $(0,\ldots,0,1,0,\ldots)$, with $1$ in the $i$th position. Our identity reduces then to the unique expansion of $m$ in base $2$.

(ii) is immediate from the definitions.

(iii) Let $[\vect a0{{t-1}}]$ be a cylinder in $Z_2^\Nbb$. By Lemma~\ref{ref8} the Lebesgue measure of its $\iota$-counterimage equals its Haar measure $\mu([\vect a0{{t-1}}])$. Since the cylinders generate the Borel $\sigma$-algebra of $Z_2^\Nbb$, the push-forward of $\lambda$ by $\iota$ is $\mu$. In particular, $\int_\Gamma u_m\,d\lambda=\int_{Z_2^\Nbb}\chi_m\,d\mu$, and (iii) amounts to the well known orthonormality of the elements of the character group.

(iv) Let $f\in C(\Gamma,\Rbb)$ and $\varepsilon>0$. Then $f\upsilon\in
C(Z_2^\Nbb,\Rbb)$ and, since the $\Rbb$-span of the characters is dense in the uniform topology, there exists $\varphi=\sum_{i=0}^s\alpha_i\chi_i$ such that $\norm{\varphi-f\upsilon}<\varepsilon$. But then $\abs{\varphi\iota(p)
-f\upsilon\iota(p)}=\abs{(\sum_{i=0}^s\alpha_iu_i)(p)-f(p)}<\varepsilon$ for every $p\in\Gamma$.
\end{proof}

\begin{remark}
Let\labell{ref16} $w_0,w_1,w_2,\ldots$ denote the classical Walsh functions on $\oi$. They are definable in a manner analogous to the one in Proposition~\ref{ref10}(i), namely by $w_m=\prod_{i=0}^{t-1}(rD^i)^{d_i}$, with the doubling map $D$ in place of the tent map $T$; the functions $rD^i$ are then the classical Rademacher functions~\cite[p.~116]{kuipersnie74}.
Let $W_t=\{w_m:2^{t-1}\le m<2^t\}$ be the set of Walsh function of level $t\ge1$. Then $W_t=U_t$ for every $t$ (we neglect the behavior at points of discontinuity). Indeed, the above identity is true for $t=1$ (since $w_1=r=u_1$) and for $t=2$ (since $w_2=rD=rT\cdot r=u_3$ and $w_3=rD\cdot r=rT=u_2$).
Now observe that:
\begin{itemize}
\item[(i)] $TD=T^2$;
\item[(ii)]
$$
U_{t+1}=\bigcup_{u\in U_t}\{uT,uT\cdot r\}.
$$
\end{itemize}
Applying the identity (ii) twice yields
$$
U_{t+1}=\bigcup_{v\in U_{t-1}}\{vT^2, vT^2\cdot rT, vT^2\cdot r, 
vT^2\cdot rT\cdot r\},
$$
and similar identities hold for $W_{t+1}$, with $D$ replacing $T$.
We thus obtain by induction
\begin{align*}
W_{t+1}&=\bigcup_{w\in W_t}\{wD,wD\cdot r\}\\
&=\bigcup_{u\in U_t}\{uD,uD\cdot r\}\\
&=\bigcup_{v\in U_{t-1}}\{vTD, vTD\cdot rD, vTD\cdot r, 
vTD\cdot rD\cdot r\}\\
&=\bigcup_{v\in U_{t-1}}\{vT^2, vT^2\cdot rT\cdot r, vT^2\cdot r,
vT^2\cdot rT\}\\
&=U_{t+1}.
\end{align*}
\end{remark}

\begin{remark}
For $n\ge3$ our $T$-Walsh functions are rather different from the usual multidimensional Walsh functions~\cite{hellekalek94}. Indeed, the latter assume constant value $\pm1$ on the cubes of a partition on $\oi^n$ defined by equations of the form $x_i=c$, with $c$ a dyadic rational. On the other hand, the $T$-Walsh function $u_m$ of level $t$ assumes constant value on the interior of the simplexes of the partition $\Bcal_t$. The equations defining the simplexes of $\Bcal_t$ have rational coefficients, but definitely not dyadic ones. Actually, pictures drawn for the $2$-dimensional case ---such as Figure~3, or the second picture on~\cite[p.~252]{panti08}, displaying $\Bcal_4$--- are rather misleading, since they lead to think that the simplexes constituting $\Bcal_t$ are always congruent to each other, and are definable by equations
of the form $\sum_{i=1}^nb_ix_i=c$, with $b_i\in\{0,1,-1\}$ and
$c\in\Zbb[1/2]$.
These facts are true for $n=1,2$, but false for $n\ge3$;
for example, one of the planes bounding $\tau_1^{11}\Gamma$ in dimension $3$ has equation $20x_1-12x_2+16x_3=15$. See also the remark after~\cite[Corollary~2.4]{panti08}.
\end{remark}

The main result in this section is the following.

\begin{theorem}
For\labell{ref12} every $p\in\Gamma$, the sequence $p,Kp,K^2p,\ldots$ is
$\lambda$-uniformly distributed.
\end{theorem}

It is well known that the $\nu$-uniform distribution of all orbits implies that a dynamical system admits only $\nu$ as an invariant measure, i.e., is uniquely ergodic. Note that this fact is usually stated for homeomorphisms (in which case the two conditions are equivalent), but the proof of the implication concerning us holds for all Borel maps; see, e.g., \cite[Theorem~I.8.2]{CornfeldFomSi82}.

By~\cite[Theorem~III.1.1]{kuipersnie74} and Proposition~\ref{ref10}(iv),
the class of $T$-Walsh functions is convergence-determining, so the validity of~\eqref{eqref2} ---with $\nu=\lambda$---
for all $T$-Walsh functions implies its validity for all continuous functions. Taking Proposition~\ref{ref10}(iii) into account, and noting
that~\eqref{eqref2} surely holds for $u_0=\one_\Gamma$, we must prove that for every $m>0$ and every $p\in\Gamma$ we have
\begin{equation}\label{eqref3}
\lim_{k\to\infty}\frac{1}{k}\sum_{i=0}^{k-1}u_m(K^ip)=0.
\end{equation}

Fix therefore $u_m$ of level $t\ge1$; then $u_m$ is constant on each $\Gammafaa$. By Lemma~\ref{ref8}, $u_m$ is constant on the topological interior of each of the $n$-dimensional simplexes in $\Bcal_t$. Let $X$ be the union of these topological interiors, and $Y=\Gamma\setminus X$ be the union of all the $(n-1)$-dimensional simplexes in $\Bcal_t$.

Choose now $p$, and let $\abar^0$ be its final orbit. Construct the sequences $\Afrak=\abar^0,\abar^1,\abar^2,\ldots$ and $\Afrak'$ as in Section~\ref{ref11}, now taking them infinite to the right only. Partition $\Afrak$ in \newword{blocks} $B_0,B_1,B_2,\ldots$ of consecutive elements: the block $B_0$ contains $\abar^0$ and all subsequent elements until and including the first one $\abar^{j_0}$ beginning with
$1^t$. The block $B_1$ now contains $\abar^{j_0+1}$ (that begins with $0^t$) and all subsequent elements till the first one $\abar^{j_1}$ beginning with $1^t$; the blocks $B_2,B_3,\ldots$ are constructed analogously. 
Let us call the last element of the block $B$ the \newword{pivot} of $B$; it is the only element of $B$ beginning with $1^t$. 

\begin{lemma}
Let\labell{ref13} $B$ be a block in $\Afrak$ with pivot $\bbar=1^t\cbar$; then the following statements hold.
\begin{itemize}
\item[(i)] If $\abar\in B$ is a final orbit, then $\bbar$ is a final orbit.
\item[(ii)] Assume that $\bbar$ is a final orbit, while $\abar\in B$ is not.
Then $\{\upsilon(\dbar):\dbar\in B\}\subset Y$.
\end{itemize}
\end{lemma}
\begin{proof}
Note first that the elements of $B$ are exactly the symbolic orbits of the form $w\cbar$, for $w$ a word of length $t$. Hence, if $w\cbar$ is a final orbit, then $1^t\cbar$ is a final orbit by Lemma~\ref{ref6}(v); this proves (i). About (ii): let $\abar=w\cbar$ be not final, $\dbar=d_0\ldots d_{t-1}\cbar\in B$, $q=\upsilon(\abar)$, $z=\upsilon(\dbar)$. Since the final orbit of $q$ is not $\abar$, while the final orbit of $T^tq$ is $\cbar$ (because $\cbar$ is a final orbit, again by Lemma~\ref{ref6}(v) applied to $\bbar$), we must have $T^iq\in\Gamma_0\cap\Gamma_1$ for some $0\le i<t$. Therefore $T^tq=T^tz$ lies in a proper face of $\Gamma$. By Lemma~\ref{ref6}(i) and the proof of~\cite[Proposition~2.2]{panti08},
$z=\tau_{d_0}\cdots\tau_{d_{t-1}}T^tz\in Y$.
\end{proof}

Rename the elements of $\Afrak'$ by writing
$$
\Afrak'=\bbar^0,\ldots,\bbar^{s_0},\bbar^{s_0+1},\ldots,\bbar^{s_1},
\bbar^{s_1+1},\ldots,\bbar^{s_2},\ldots
$$
where $\bbar^{s_i}$ is the $i$th surviving pivot; since $\abar^0$ is a final orbit, we have $\bbar^0=\abar^0$ and $\bbar^{s_0}$ is the pivot of $B_0$,
by Lemma~\ref{ref13}(i).
The \newword{reduced blocks} $B'_i$ are defined by
$B'_0=\{\bbar^0,\ldots,\bbar^{s_0}\}$ and
$B'_{i+1}=\{\bbar^{s_i+1},\ldots,\bbar^{s_{i+1}}\}$. The elements of $B'_i$
are exactly the final orbits in the block $B\subset\Afrak$ to which $\bbar^{s_i}$ belongs.
The cardinality of $B'_i$ varies from $1$ (if only the pivot survives) to $2^t$. Applying Theorem~\ref{ref2}, let $p^i=K^ip$ be the point whose final orbit is $\bbar^i$,
and let $S_i=\sum\{u_m(p^j):\bbar^j\in B'_i\}$. Since the length of the reduced blocks is bounded, and $u_m$ is bounded too, we can establish~\eqref{eqref3} by computing the limit along the surviving pivots, i.e., by establishing
\begin{equation}\label{eqref4}
\lim_{k\to\infty}\frac{1}{s_k+1}\sum_{i=0}^kS_i=0.
\end{equation}
Consider the reduced block $B'_i$. If $B'_i$ contains $2^t$ elements, then one sees easily that $S_i=0$. By Lemma~\ref{ref13}(ii), if $B'_i$ contains less than $2^t$ elements, then $\{p^j:\bbar^j\in B'_i\}\subset Y$. Therefore~\eqref{eqref4} follows immediately from the following claim.

\noindent\emph{Claim.} The set $\{j\in\Nbb:p^j\in Y\}$ has density $0$.

\noindent\emph{Proof of Claim.} Let us backtrack to the unpruned sequence $\Afrak$, and let $q^i=\upsilon(\abar^i)$. The sequence $q^0,q^1,q^2,\ldots$ may have repetitions, but is surely $\lambda$-uniformly distributed. Indeed, every continuous function $f:\Gamma\to\Rbb$ gives rise to a continuous function $f\upsilon$ on $\Zbb_2$. One then uses the fact that $\Afrak$ is
$\mu$-uniformly distributed~\cite[Theorem~IV.4.2]{kuipersnie74}, and that
$\lambda$ is the push-forward of $\mu$ by $\upsilon$, as noted in Corollary~\ref{ref17}.

Choose now $\varepsilon>0$. By the above, we can find an index $i_0$ such that for each $i_1\ge i_0$ we have
$$
\frac{\sharp\{0\le i<i_1:q^i\in Y\}}{i_1}<\frac{\varepsilon}{2^{n(n-1)/2}},
$$
where $\sharp$ denotes cardinality. Without loss of generality $\abar^{i_0}$ is a surviving pivot, say $\abar^{i_0}=\bbar^{j_0}$. Let $j_1\ge j_0$, so that $\bbar^{j_1}=\abar^{i_1}$ for a certain $i_1\ge i_0$.
By Theorem~\ref{ref1}, every point in $q^0,q^1,q^2,\ldots$ repeats at most $2^{n(n-1)/2}$ times. Therefore
$j_1\ge i_1/2^{n(n-1)/2}$, and we obtain
\begin{multline*}
\frac{\sharp\{0\le j<j_1:p^j\in Y\}}{j_1}\le
2^{n(n-1)/2}\frac{\sharp\{0\le j<j_1:p^j\in Y\}}{i_1}\\
\le 2^{n(n-1)/2}\frac{\sharp\{0\le i<i_1:q^i\in Y\}}{i_1}<\varepsilon.
\end{multline*}
This proves our claim, and concludes the proof of Theorem~\ref{ref12}.

\section{An arithmetical conjugate}

In this final section we discuss a conjugate of our map $K$ that has some arithmetical significance. Since the facts we are presenting derive in a rather formal fashion from the results proved in the previous sections, we will be somehow brief.

Recall the matrices $C_0$ and $C_1$ introduced in Section~\ref{ref7}. Define a map $\psi_0:\Gamma\to\Gamma$ as follows: if $p=(\vect\alpha1n)\in\Gamma$, then $\psi_0p$ is the unique point $(\vect\beta1n)$ such that $(\beta_1\cdots\beta_m\,1)^{tr}$ is proportional to $C_0(\alpha_1\cdots\alpha_m\,1)^{tr}$.
Define analogously $\psi_1$ in terms of $C_1$.
For every $t\ge1$ and every $t$uple $(\vect a0{{t-1}})\in\ooii^t$ the image 
$\psi_{a_0}\cdots\psi_{a_{t-1}}\Gamma$ is a simplex, and the set of all faces of these $2^t$ $n$-dimensional simplexes form a simplicial complex $\Fcal_t$ supported on $\Gamma$. The complexes $\Fcal_t$ and $\Bcal_t$ are combinatorially isomorphic, and there exists a unique
orientation-preserving homeomorphism
$\Phi:\Gamma\to\Gamma$ that restricts to homeomorphisms between $\psi_{a_0}\cdots\psi_{a_{t-1}}\Gamma$ and 
$\tau_{a_0}\cdots\tau_{a_{t-1}}\Gamma$, for each $t$ and $\vect a0{{t-1}}$. 
See~\cite{panti08} for the above results; the map $\Phi$ can be seen as an $n$-dimensional generalization of the Minkowski question mark function~\cite{kinney60}, \cite{viaderparjabi98}.

Since $K$ was defined via a combinatorial property (namely, the existence of the partitions in Proposition~\ref{ref5}), it is no surprise that the conjugate $\Phi^{-1} K\Phi$ is definable via an analogous combinatorial construction. Namely, we define a bijection $E:\Gamma\to\Gamma$ by setting
$E=\psi_0^k\psi_1(\psi_1^k\psi_0)^{-1}$ on $\psi_1^k\psi_0\Gammao$.
We also set $Ev'_{-1}=v_0$, where $v'_{-1}$ is the only element in $\bigcap\{\psi_1^k\Gamma:k\ge0\}$. By~\cite[Proposition~3.1]{panti08} the following diagram commutes
\begin{equation*}
\begin{xy}
\xymatrix{
\psi_1^k\psi_0\Gammao \ar[d]_\Phi &
\Gammao \ar[l]_{\psi_1^k\psi_0} \ar[d]^\Phi \ar[r]^{\psi_0^k\psi_1} &
\psi_0^k\psi_1\Gammao \ar[d]^\Phi \\
\tau_1^k\tau_0\Gammao &
\Gammao \ar[l]^{\tau_1^k\tau_0} \ar[r]_{\tau_0^k\tau_1} &
\tau_0^k\psi_1\Gammao
}
\end{xy}
\end{equation*}
Since $\Phi v'_{-1}=v_{-1}$ and $\Phi v_0=v_0$, we have $E=\Phi^{-1}K\Phi$, as expected.

The bijection $E$ is piecewise-fractional with integer coefficients: indeed, let $(e^i_1\cdots e^i_{n+1})$ be the $i$th row of $C_0^kC_1C_0^{-1}C_1^{-k}$. Then on $\psi_1^k\psi_0\Gammao$ the $i$th component $E^i$ of $E$ (i.e., $E$ followed by the projection on the $i$th coordinate) has the form
$$
E^i(\vect\alpha1n)=
\frac{e^i_1\alpha_1+\cdots+e^i_n\alpha_n+e^i_{n+1}}{e^{n+1}_1\alpha_1+\cdots+e^{n+1}_n\alpha_n+e^{n+1}_{n+1}}.
$$
In the $1$-dimensional case a conjugate of the classical Kakutani-von Neumann map via the Minkowski function was introduced in~\cite[Theorem~2.3]{bonannoisola09}.

\begin{theorem}
The\labell{ref14} homeomorphism $E$ is minimal and uniquely ergodic, with the Minkowski measure $\Phi^*\lambda$ as its unique invariant probability. All points of $\Gamma$ have a $\Phi^*\lambda$-uniformly distributed $E$-orbit. The orbit of $v_0$ constitutes an enumeration without repetitions of all points in $\Gamma$ having rational coordinates.
\end{theorem}

\begin{proof}
By definition $\Phi^*\lambda$ is the pullback of $\lambda$ via $\Phi$, i.e., $(\Phi^*\lambda)(A)=\lambda(\Phi A)$ for every Borel subset $A$ of $\Gamma$. All statements are immediate from our previous results, upon noting that by~\cite[Theorem~3.5]{panti08} the set of rational points in $\Gamma$ is mapped bijectively by $\Phi$ to the set of dyadic points.
\end{proof}

Basic ergodic theory yields that the Minkowski and the Lebesgue measures are mutually singular~\cite[p.~262]{panti08}.
The enumeration of all rational points given by Theorem~\ref{ref14} gives a nice and effective representation of the mass distribution determined by the Minkowski measure, so we conclude this paper by drawing the first 6000 points in the $E$-orbit of $v_0$.
\begin{figure}[H]
\caption{the set $\{E^tv_0:0\le t<6000\}$.}
\begin{center}
\includegraphics[height=7cm]{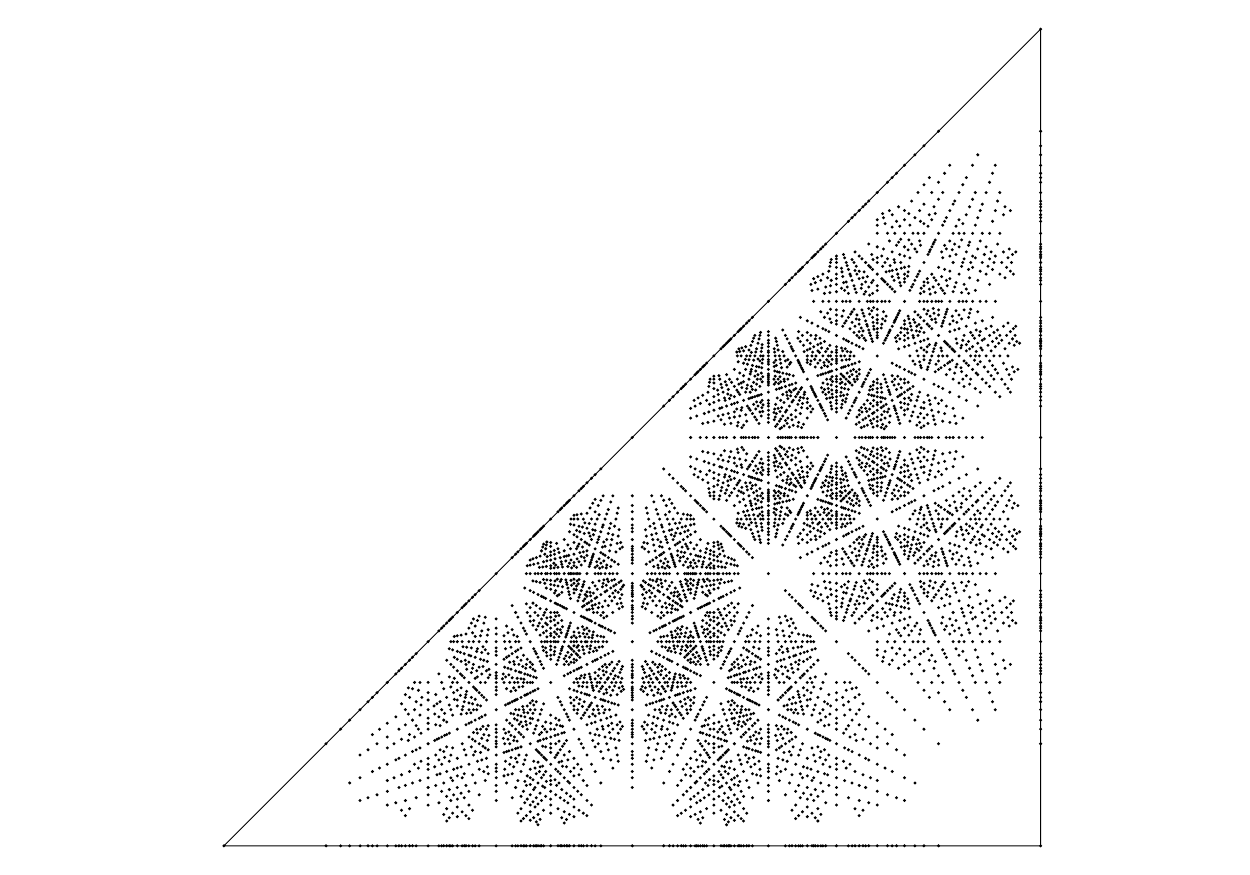}
\end{center}
\end{figure}

\end{document}